\documentclass[10pt,twoside]{siamltex}
\usepackage{amsfonts,epsfig,amsmath,amssymb,graphicx,geometry}

\def\IC{\hbox{\rm C\kern-.43em
       \vrule depth 0ex height 1.4ex width .05em\kern.41em}}
\def\IQ{\hbox{\rm Q\kern-.43em
       \vrule depth 0ex height 1.4ex width .05em\kern.41em}}

\def\qed{\hfill\vbox{\hrule width 6 pt
\hbox{\vrule height 6 pt width 6 pt}}}

\begin{document}

%  Leave these commented lines here
% \input{elaheader-volx-xx.tex}
% \setcounter{page}{1}

% \renewcommand{\thefootnote}{\fnsymbol{footnote}}
% \renewcommand{\thefootnote}{\arabic{footnote}}
% \renewcommand{\theequation}{\thesection.\arabic{equation}}

\bibliographystyle{plain}
\title{Factorization of permutations}
% Leave blank; editors will write the exact dates above

\author{
Zejun Huang\thanks{College of Mathematics and Econometrics, Hunan University, Changsha, Hunan 410082, China.
(Email: mathzejun@gmail.com)   }
% Remember to put \and between any two authors
\and
Chi-Kwong Li\thanks{Department of Mathematics, College of William and Mary, Williamsburg, VA 23187,
USA; Department of Mathematics, University of Hong Kong, Pokfulam, Hong Kong. (Email:
ckli@math.wm.edu)  }
\and
Sharon H. Li\thanks{Department of Computer Science, The Johns Hopkins University, Baltimore, MD 21218,
USA. Current affiliation: Microsoft Corporation, Redmond, WA 98052.
(sharonli@jhu.edu)}
\and
Nung-Sing Sze\thanks{Department of Applied Mathematics, The Hong Kong Polytechnic University, Hung
Hom, Hong Kong. (raymond.sze@polyu.edu.hk)  }}

% Authors and running title to go on top of each page
\pagestyle{myheadings}
\markboth{Z. Huang, C.K. Li, S.H. Li, and N.S. Sze}{ Factorization of permutations}
\maketitle

\begin{abstract}
We consider the problem of factoring permutations as a product of
special types of transpositions, namely, those transpositions
involving two positions with bounded distances. In particular, we
investigate the minimum number, $\delta$,
such that every permutation can be factored into no more than
$\delta$ special transpositions. This study is related to
sorting algorithms, Cayley graphs, and genomics.
\end{abstract}

%%%%%%%%%%%%%%%%%%%%%%%%%%%%%%%%%%%%%%%%%%%%%%%%%%%%%%%%%%%%%

\section{Introduction}
A basic problem in computer science concerns the sorting of a list of elements in a random order
to a specific order. For example, the bubble sort algorithm can be used to restore the order
of a list of numbers, say, $[i_1, \dots, i_n]$, an arrangement of the numbers $[1, \dots, n]$,
by swapping adjacent elements to restore the list to its natural (ascending) order.

Mathematically, we identify the list $\sigma = [i_1, \dots, i_n]$ as a permutation in $S_n$,
the symmetric group of degree $n$, such that $\sigma(j) = i_j$ for $j = 1, \dots, n$.
Define the number of inversions of the
permutation $\sigma = [i_1, \dots, i_n]$ as
the sum of the numbers inv$(j)$, where inv$(j)$ is the number
of integers smaller than $j$ lying on the right side
of $j$ in $[i_1, \dots, i_n]$. For example, the number of
inversions in $[3,2,4,5,1]$ is $2+1+1+1+0 = 5$.

Applying bubble sort to the permutation $\sigma = [i_1, \dots, i_n]$ corresponds to restoring
$\sigma$ to the identity permutation $[1, \dots, n]$ by exchanging
two adjacent numbers in each step.
It is not hard to see that the minimum number of steps needed
is the number of inversions in $\sigma = [i_1, \dots, i_n]$.
In fact, switching two adjacent numbers of a permutation will
increase or decrease the number of inversions by one.
So, for $\sigma\in S_n$, if in each step one swaps two adjacent numbers that are in
the wrong order, i.e., so that the left one is larger than the right one, then
one will get identity permutation after $k$ steps, where
$k$ is the number of inversions of $\sigma$.
Hence, the worst scenario is when the permutation
$[n,n-1, \dots, 1]\in S_n$, which has the maximum number of inversions:
$(n-1) + \cdots + 1= n(n-1)/2$.
There are many other efficient sorting algorithms.
We refer the readers to \cite{Knuth}  for more details.

In this note, we consider the problem of finding the minimum number
of steps needed to convert a permutation to the identity permutation
if one is allowed to switch the numbers in the $i$th and $j$th positions
as long as $|i-j| \le m$, for some $m \in \{1, \dots, n - 1\}$.
Clearly, the bubble sort algorithm is the case when $m = 1$.

Let $G_m$ be the set of transpositions $(i,j)$ in $S_n$ with $j-i \le m$.
We investigate the minimum number $\delta(n,m)$ such that every
permutation can be factored into no more than $\delta(n,m)$ transpositions in $G_m$.
For $(n-1)/2 \le m$, we give a formula for $\delta(n,m)$, and
characterize all permutations in $S_n$ requiring  $\delta(n,m)$
transpositions in its factorization (see Theorem \ref{thm2.6}). For
$1 < m < (n-1)/2$, we obtain an upper bound for $\delta(n,m)$ (see
Section 2.2).

Note that the result on bubble sort can be formulated in terms of
the Cayley graph of $S_n$, constructed as follows:
Represent every permutation as a vertex, and
connect two vertices $\sigma_1$ and $\sigma_2$
if $\sigma_2 = \tau \sigma_1$ for a permutation $\tau = (i,i+1)$
which exchanges two adjacent numbers at the $i$th and $(i+1)$th positions
for some $i \in \{1, \dots, n-1\}$. Then the bubble sort algorithm
amounts to moving a permutation $[i_1, \dots, i_n]$ to $[1, \dots, n]$
in the Cayley graph most efficiently (using the minimum number of steps).
Moreover, $n(n-1)/2$  is the maximum distance from $[1, \dots, n]$ to another
permutation $[i_1, \dots, i_n]$, which is $[n,\dots,1]$.
One easily shows that the value $n(n-1)/2$ is actually the maximum
distance between any two vertices, and is known as the diameter of the Cayley graph.
It also indicates that a permutation is a product of no more than
$n(n-1)/2$ transpositions of the form $(i,i+1)$.

Our study concerns the Cayley graph of using elements in $S_n$ as vertices
so that two vertices $\sigma_1$ and $\sigma_2$ are connected
if $\sigma_2 = \tau \sigma_1$ for a permutation $\tau \in G_m$.

The study of Cayley graphs of $S_n$ has a long history; see \cite{JD,VF, GF,IR,L,M,xiao} and
their references. Note that in some of these papers, the authors study the minimum number of
transpositions needed in the factorization. It is in essence studying the
diameter of the underlying Cayley graphs.

It is interesting to note that the study is related to other topics such as genome rearrangement.
In nature, some species have similar genetic make up and differ only in the
order of their genes. Finding the shortest rearrangement path
between two related bacteria or viruses is useful in drug discovery and vaccine
development. The study is also useful in the study of mutations.
In fact, a slight change of the genetic sequence
may have significant effect, and it is more likely to see a change (permutation)
of the positions of the nucleotides close to each other in the genetic sequences.
That is why there is keen interest in studying
such permutations in genetic sequences; for example, see
\cite{Feng,P,Yue} and their references.

\section{Main results}

Following Section 1, for $m \in \{1, \dots, n-1\}$, let  $G_m$ be the set of
transpositions $(i,j)$ in $S_n$ with $j-i \le m$.
Then $G_m$ generates $S_n$, i.e., every permutation in $S_n$ is a product of
transpositions in $G_m$.

We are interested in finding the smallest number of transpositions
in $G_m$ needed to convert a given permutation to the identity,
and those permutations which require the maximum number of
transpositions to do the reduction.

Consider the Cayley graph $\Gamma_{n,m}$  so
that the vertices are elements in $S_n$, and two vertices $\sigma_1$
and $\sigma_2$ if $\sigma_1\sigma_2^{-1} \in G_m$.
Denote by {\bf 1} the identity permutation $[1,\dots,n]$.
We are interested in the shortest path connecting the identity permutation
${\bf 1}$ to a given permutation $\sigma$. The length of this path is denoted by
$d({\bf 1}, \sigma,m)$. Also, we are interested in
the permutation $\sigma^*$ with a maximum distance to the identity ${\bf 1}$.
Clearly, the maximum distance will be the same as the maximum
distance between any two vertices in the Cayley graph, and we will denote this quantity
by $\delta(n,m)$, as introduced in Section 1.
By the discussion in Section 1, we have $\delta(n,1) = n(n-1)/2$.

Let $\sigma = (j_1,j_2,\ldots,j_r) \in S_n$ be a cycle of length $r$, i.e.,
the permutation $\sigma\in S_n$ such that
$\sigma(j_1)=j_2,\ldots,\sigma(j_{r-1})=j_r,\sigma(j_r)=j_1$
and $\sigma(j) = j$ for other $j$.
Also denote by $\lceil x \rceil$ the smallest integer greater than or equal to  $x$
and $\lfloor x \rfloor$ the largest integer less than or equal to $x$ respectively.

%\subsection{The case when $m = n-1$}

The following lemma will be used frequently in our discussion.

\begin{lemma}\label{le1}
Let $C_1=(i_1,i_2,\cdots,i_k)$ and $C_2=(j_1,j_2,\cdots,j_s)$ be two disjoint cycles in $S_n$.
Suppose $\tau_1=(i_\alpha,i_\beta)$ and $ \tau_2=(i_u,j_v)$ with
$1\leq\alpha<\beta\leq k$, $1\leq u\leq k$ and $1\leq v\leq s$. Then
$$
\tau_1C_1 =(i_\alpha,i_\beta)(i_1,i_2,\cdots,i_k)
      = (i_1,i_2,\cdots,i_{\alpha-1},i_\beta,i_{\beta+1},\cdots,i_k)
			(i_\alpha,i_{\alpha+1},\cdots,i_{\beta-1})$$
is the product of two disjoint cycles, and
\begin{eqnarray*}
\tau_2C_1C_2&=&(i_u,j_v)(i_1,i_2,\cdots,i_k)(j_1,j_2,\cdots,j_s)\\
        &=&(i_1,i_2,\cdots,i_{u-1},j_v,j_{v+1},\cdots,j_s,j_1,j_2,\cdots,j_{v-1},i_u,i_{u+1},\cdots,i_k)
\end{eqnarray*}
is a cycle.
\end{lemma}

\vskip 8pt
The result for $m = 1$ was discussed in Section 1.
The other extreme is when $m = n-1$, i.e., $G_m$ is the set of all transpositions.
We have the following known result, see \cite{JD, M}.
Here we give a short proof of it for completeness.

\begin{proposition} \label{2.2} Suppose $\sigma \in S_n$. Then
$d({\bf 1}, \sigma, n-1) = n-r$, where $r$ is the number of cycles in the disjoint
cycle representation of $\sigma$ under the convention that each fixed point
is counted as a  1-cycle. Thus,
$\delta(n, n-1) = n-1$ is attained at a  $n$-cycle.
\end{proposition}

\it Proof. \rm
Writing $\sigma = \tau_1 \cdots  \tau_k\in S_n$ for a minimum number of
transpositions $\tau_1, \dots, \tau_k$ is the same as
finding the minimum number of transpositions $\tau_1 \cdots \tau_k$
such that $ \tau_k \cdots \tau_1\sigma = {\bf 1}$. Using Lemma \ref{le1},
to convert $\sigma$ to the product of $n$ disjoint cycles, i.e.,
back to the identity permutation $[1, \dots, n]$ most efficiently by
composing $\sigma$ with transpositions, the most efficient way is to
choose transposition $(i,j)$ in each step such that $i$ and $j$ lie in
the same cycle. So, using $n-r$ transpositions $\tau_1, \dots, \tau_{n-r}$
to convert $\sigma$ to the identity permutation will be the most efficient
scheme. \qed

%Denote by $\lceil x \rceil$ the smallest integer greater than or equal to  $x$
%and $\lfloor x \rfloor$ the largest integer less than or equal to $x$ respectively.

\subsection{The case when $m \ge  (n - 1)/2$}

Given $\sigma\in S_n$, denote by $K_1(\sigma)$ the set of transpositions in
$G_{n-1}$ splitting a cycle of $\sigma$ into two and by $K_2(\sigma)$
the set of transpositions in $G_{n-1}$ jointing two cycles of  $\sigma$  into one. We call
$K_1(\sigma)$ and $K_2(\sigma)$ type one and type two transpositions, respectively.
By Lemma \ref{le1} we have

\begin{proposition}\label{le7}
Let $\sigma\in S_n$, $\tau_1\in K_1(\sigma)$ and $\tau_2\in K_2(\sigma)$. Then
$$d({\bf 1},\tau_1\sigma,n-1)=d({\bf 1}, \sigma,n-1)-1\quad
\hbox{ and } \quad
  d({\bf 1},\tau_2\sigma,n-1)=d({\bf 1}, \sigma,n-1)+1.$$
 \end{proposition}

For any cycle $C = (i_1,\dots,i_p)$ in $S_n$,
we say that $C$ is in the set $L_m$ if
for each term $i_u$ in $C$, there is a term $i_v$ such that
$|i_u - i_v| > m$. On the other hand,
a cycle $C\notin L_m$ if
there is some $1\le r \le p$ such that
$|i_t - i_r| \le m$ for all $t= 1,\dots,p$.

\begin{lemma}\label{2.9}
Let $n$ and $m$ be positive integers and $5 \le n \le 2m+1$.
Suppose $C = (i_1,i_2,\dots, i_p)$ is a length $p$ cycle in $S_n$.
Then one of the following holds.
\begin{enumerate}

\item[\rm (a)] If $C \in L_m$,
then $C$ can be written as a product of $p + 1$ transpositions in $G_{m}$.

\item[\rm (b)] If $C \notin L_m$,
then $C$ can be written as a product of $p-1$ transpositions in $G_m$.
\end{enumerate}
Furthermore, suppose $C_i = (i_1,\dots,i_p)$ and $C_j = (j_1,\dots,j_q)$ are two cycles in $L_m$.
If there are $1\le r \le p$ and $1\le s \le q$ such that
\begin{eqnarray}\label{rs}
|i_t - j_s| \le m \quad\hbox{for all}\quad t = 1,\dots,p
\quad\hbox{and}\quad
|j_u-i_r| \le m \quad\hbox{for all}\quad  u = 1,\dots,q,
\end{eqnarray}
then $C_i C_j$ can be written as a product of $p+q$ transpositions in $G_m$.
\end{lemma}

\it Proof. \rm
Suppose $C \in L_m$.
Notice that $|i_t  - (m+1)| \le m$ for all $t$.
Then one can write
$$C = (m+1,i_p)\, (m+1,i_{p-1})\, (m+1,i_{p-2})\, \cdots\, (m+1,i_2)\, (m+1,i_1)\, (m+1,i_p).$$
Then the result (a) holds.

Suppose now $C \notin L_m$. That is,
there is some $1\le r \le p$ such that
$|i_t - i_r| \le m$ for all $t= 1,\dots,p$.
Without loss of generality,
we may assume that $r = p$. Then $C$ can be written as
$$C = (i_p,i_{p-1})\, (i_p,i_{p-2})\, (i_p,i_{p-3})\, \cdots\,  (i_p,i_{2})\, (i_p,i_{1}).$$
Thus, the result (a) holds.

Suppose now $C_i$ and $C_j$ are disjoint cycles and satisfying (\ref{rs}).
we may assume that $r = p$ and $s = q$ in (\ref{rs}). Then $C_i C_j$ can be written as
$$C_iC_j =  (j_q,i_p)\, (j_q,i_{p-1})\, (j_q,i_{p-2})\, \cdots\,  (j_q,i_{2})\, (j_q,i_{1})\,
(i_p,j_{q-1})\, (i_p, j_{q-2})\,\cdots (i_p,j_1)\, (i_p,j_q).$$
Thus, the result follows.
\qed

\begin{lemma}\label{2.10}
Let $n$ and $m$ be positive integers and $5 \le n \le 2m+1$.
Suppose $\sigma \in S_n$ has a disjoint cycle decomposition
$C_1\cdots C_r$ under the convention that each fixed point is counted as a $1$-cycle.
If $s$ of the cycles $C_i$ belongs to $L_m$
and among them, $t$ disjoint pairs of $C_i$ and $C_j$ satisfy condition (\ref{rs}) in Lemma \ref{2.9},
then $d({\bf 1},\sigma,m) \le n -r + 2s - 2t$.
\end{lemma}

\it Proof. \rm
Suppose $\sigma \in S_n$ has a disjoint cycle decomposition $C_1\cdots C_r$
such that $C_j \in L_m$ for $j = 1,\dots, s$,
$C_j \notin L_m$ for $j = r-s+1,\dots,r$.
Further,  the cycles $C_{2k-1}$ and $C_{2k}$
satisfy condition (\ref{rs}) for $k = 1,\dots, t$.

Assume that $C_j$ has length $\ell_j$ for $j=1,\ldots,r$.
Then by Lemma \ref{2.9},
\begin{eqnarray*}
d({\bf 1},\sigma,m)
&\le&   \sum_{j=1}^{2t} \ell_j  + \sum_{j=2t+1}^{s} (\ell_j + 1) + \sum_{j=s+1}^r (\ell_j - 1)\\
&=& \left(\sum_{j=1}^r \ell_j \right) +(s-2t) - (r-s)
= n -r + 2s-2t.
\end{eqnarray*}
\vskip -.5in   \qed

\vskip .3in

\begin{theorem}
\label{thm2.6}
Let $n$ and $m$ be positive integers and $5 \le n \le 2m+1$. Then
$$\delta(n,m) = n+d- 1 \quad \hbox{
with } \quad d = \left\lfloor \frac{n-m}{2} \right\rfloor.$$
A permutation $\sigma \in S_n$
attains $d({\bf 1},\sigma,m) = n+ d - 1$
if and only if one of the following holds.
\begin{enumerate}
\item[\rm (a)] $(n-m)$ is even and $\sigma$ is a product of $d+1$ disjoint cycles of the form
$$(i_1,j_1) (i_2,j_2) \cdots (i_d,j_d) (k_1,\dots,k_{n-2d}),$$
where
$\{i_1,\dots,i_d\} = \{1,\dots,d\}$,
$\{j_1,\dots,j_d\} = \{n-d+1,\dots,n\}$,
and
$\{k_1,\dots,k_{n-2d}\} = \{d+1,\dots,n-d\}$.

\medskip
\item[\rm (b)] $(n-m)$ is odd and $\sigma$ is a product of $d+1$ disjoint cycles of the form

\medskip
\item[\rm (b.1)] \ $ (i_1,j_1) (i_2,j_2) \cdots (i_{d-1},j_{d-1}) (i_d,j_d)
(k_1,\dots,k_{n-2d}),$

\medskip
where
$\{i_1,\dots,i_d\} \subseteq \{1,\dots,d+1\}$,
$\{j_1,\dots,j_d\} \subseteq \{n-d,\dots,n\}$,
and \linebreak
$\{k_1,\dots,k_{n-2d}\} \subseteq \{d+1,\dots,n-d\}$
such that
$\{d+1,n-d\} \cap \{k_1,\dots,k_{n-2d}\} \ne \emptyset$.

\medskip
\item[\rm (b.2)] \ $(i_1,j_1) (i_2,j_2) \cdots (i_{d-1},j_{d-1}) (i_d,j_d,i_{d+1})
(k_1,\dots,k_{n-2d-1}),$

\medskip
where
$\{i_1,\dots,i_{d+1}\} = \{1,\dots,d+1\}$,
$\{j_1,\dots,j_d\} = \{n-d+1,\dots,n\}$, and \linebreak
$\{k_1,\dots,k_{n-2d-1}\} = \{d+2,\dots,n-d\}$.

\medskip
\item[\rm (b.3)] \ $(i_1,j_1) (i_2,j_2) \cdots (i_{d-1},j_{d-1}) (i_d,j_d,j_{d+1})
(k_1,\dots,k_{n-2d-1}),$

\medskip
where
$\{i_1,\dots,i_{d}\} = \{1,\dots,d\}$,
$\{j_1,\dots,j_{d+1}\} = \{n-d,\dots,n\}$, and \linebreak
$\{k_1,\dots,k_{n-2d-1}\} = \{d+1,\dots,n-d-1\}$.

\medskip
\item[\rm (b.4)] \
$(i_1,j_1) (i_2,j_2) \cdots (i_{d-1},j_{d-1}) (i_d,j_d,i_{d+1},j_{d+1})
(k_1,\dots,k_{n-2d-2}),$

\medskip
where
$\{i_1,\dots,i_{d+1}\} = \{1,\dots,d+1\}$,
$\{j_1,\dots,j_{d+1}\} = \{n-d,\dots,n\}$, and \linebreak
$\{k_1,\dots,k_{n-2d-2}\} = \{d+2,\dots,n-d-1\}$
such that
 $\{d+1,n-d\} \cap \{i_d,j_d,i_{d+1},j_{d+1}\}\ne \emptyset$.

 \end{enumerate}
\end{theorem}

 \it Proof. \rm
Suppose $\sigma$ has a disjoint cycles decomposition $C_1\cdots C_r$.
Assume $C_1,\dots,C_s$ are disjoint cycles in $L_m$
while $C_{s+1},\dots,C_r$ are not in $L_{m}$.
Notice that $|i- \lceil\frac{n+1}{2}\rceil|\le m$ for all $i=1,\ldots,n$. It follows that the
cycle containing $\lceil\frac{n+1}{2}\rceil$ is not in $L_m$ and $r>s$.
For $j=1,\dots,s$, let $u_j$ and $v_j$ be the smallest term and largest term
of the cycle $C_j$ respectively. Since $C_i \in L_m$,
 $|v_i - u_i| > m$, and
$$\{u_1,\dots,u_s\} \subseteq \{1,\dots,n-m-1\}
\quad\hbox{and}\quad \{v_1,\dots,v_s\} \subseteq \{m+2,\dots,n\}.$$
Note that if $|v_j - u_i| \le m$ then $C_i$ and $C_j$ satisfy (\ref{rs}).
Moreover, since $n-m-1<m+2$, we have $$\min_{1\le i\le s}v_i>\max_{1\le i\le s}u_i.$$

Now assume there are $t$ disjoint pairs of cycles satisfying (\ref{rs}).
Without loss of generality, we assume that no pair satisfying (\ref{rs}) can be found among
the cycles $C_1,\dots,C_{\hat s}$ with $\hat s = s - 2t$.
By the claim, we must have $|v_j - u_i| > m$ for all $1\le i , j \le \hat s$. Then
$$\hat s \le \max_{1\le i \le \hat s} u_i \le \min_{1\le j \le \hat s} v_j  - (m +1)
\le (n-\hat s+1) - (m+1).$$
Thus, $2\hat s \le n-m$ and hence $\hat s \le \left\lfloor \frac{n-m}{2} \right\rfloor = d$.
By  Lemma \ref{2.10} and the fact that $r > s$,
\begin{equation}\label{e16}
d({\bf 1}, \sigma,  m)\leq n -r + 2s -2t = n + \hat s - (r-s) \le n + d - 1.
\end{equation}
Furthermore,  equality holds only if $\hat s = d$ and $r-s = 1$.

Assume now that $\sigma$ attains the upper bound.
Then $\hat s = d$.
As any two cycles of $C_1,\dots,C_{\hat s}$ do not satisfy (\ref{rs}),
we must have
\begin{equation}\label{eq20151}
\max_{1\le i \le d} u_i \le d+1
\quad\hbox{and}\quad
\min_{1\le j \le d} v_j \ge n-d.
\end{equation}
Furthermore, at most one of the inequalities is actually an equality if $n-m$ is odd while
both two inequalities are strictly inequalities if $n-m$ is even.
If $t > 0$, the union of the two sets
$$\{u_{d + 1},\dots,u_{d + 2t}\} \cap \{1,\dots,d\}
\quad\hbox{and}\quad
\{v_{d+1}, \dots, v_{d + 2t}\} \cap \{n-d+1,\dots,n\}$$
contains at most one element. Therefore,
there is $k$ such that $u_{d+k} > d$ and $v_{d+k} < n-d+1$.
But then $$ v_{d+k} - u_{d+k} < n-d-(d+1)=n-2d-1 \le m,$$
which contradicts that $C_{d+k} \in L_m$.
Hence, $t = 0$.
Thus, $\hat s = s = d$ and $\sigma$ has a $d+1$ disjoint cycle decomposition $C_1\cdots C_{d+1}$
with $\{C_1,\dots,C_d\}\subseteq L_m$ and $C_{d+1} \notin L_m$.

If $n-m$ is even, then $2d = n-m$ and by (\ref{eq20151}) we have
$$\max_{1\le i \le d} u_i  = d \quad\hbox{and}\quad \min_{1\le j \le d} v_j  = n-d+1.$$
It follows that
\begin{eqnarray}\label{uv0}
\{u_1,\dots,u_d\} =  \{1,\dots,d\}
\quad\hbox{and}\quad
\{v_1,\dots,v_d\} =  \{n-d+1,\dots,n\}.
\end{eqnarray}
Thus, each $C_i$ has exactly one element in $\{1,\dots,d\}$
and one element in $\{n-d+1,\dots,n\}$.

Suppose any of $C_1,\dots,C_d$ has length greater than $2$.
Without loss of generality, assume $C_1 = (i_1,\dots,i_p)$ has length   $p > 2$,
with $d+1 \le i_p \le n-d$. By symmetry, let us first assume that $d+1\le i_p \le m$.
For the case for $m < i_p \le n-d$, one can obtain the same conclusion by a similar argument.
Let $i_\ell$ be the only element of $C_1$ that lies in $\{n-d+1,\dots,n\}$.
Then one can see that
$$|i_t - i_p| \le m\quad\hbox{for} \quad t = 1,\dots,p,\ t \ne \ell.$$
Since $C_1\in L_m$, we have $i_l-i_p>m$, i.e., $i_l>i_p+m\ge d+1+m=n-d+1$. By (\ref{uv0}),
there is another cycle, say $C_2$, with $v_2=n-d+1$ such that $|v_2 - i_p | \le m$.
Let $C_2 = (j_1,\dots,j_q)$. Then
$$|j_t - i_p| \le m \quad\hbox{for}\quad t = 1,\dots,q.$$
We assume that $j_q = v_2$. Then $C_1 C_2$ can be written as
$$C_1 C_2 = (i_p,i_{p-1}) \cdots (i_p,i_{\ell+1}) (i_p,j_q) (j_q,i_\ell) (i_p,j_{q-1}) \cdots (i_p,j_1)
(i_p,j_q) (i_p,i_{\ell-1}) \cdots (i_p,i_1),
$$
which is a product of $p+q$ transpositions in $G_{m}$.
By Lemma \ref{2.9}, $C_3 \cdots C_{d+1}$ can be written
as a product of $n-(p+q) + (d-3)$ transpositions in $G_m$.
Thus, $\sigma$ is a product of $n + d -3$ transpositions,
which contradicts that $\sigma$ attains the upper bound.
Therefore, all $C_i$ have length $2$ and the case (a) holds by (\ref{uv0}).

Now if $n - m$ is odd, then $2d = m-n-1$. By (\ref{eq20151}) we have  either
%\begin{eqnarray}\label{2.3}
$$\hbox{\rm (i) }
\max_{1\le i \le d} u_i  = d \quad\hbox{and}\quad \min_{1\le j \le d} v_j  \ge n-d$$
 or $$\quad\hbox{\rm (ii) }\max_{1\le i \le d} u_i  \le d+1
 \quad\hbox{and}\quad \min_{1\le j \le d} v_j  = n-d+1.
$$%\end{eqnarray}
Then either
\begin{eqnarray}\label{uv1}
\hbox{\rm (i)}\quad \{u_1,\dots,u_d\} =  \{1,\dots,d\}
\quad\hbox{and}\quad
\{v_1,\dots,v_d\} \subseteq  \{n-d,\dots,n\}
\end{eqnarray}
or
\begin{eqnarray}\label{uv2}
\hbox{\rm (ii)}\quad \{u_1,\dots,u_d\} \subseteq  \{1,\dots,d+1\}
\quad\hbox{and}\quad
\{v_1,\dots,v_d\} =  \{n-d+1,\dots,n\}.
\end{eqnarray}
Suppose any of $C_i$, $1\le i\le d$, contains an element in $\{d+2,\dots,n-d-1\}$.
Without loss of generality, assume $C_1 = (i_1,\dots,i_p)$ is the cycle and
$d+2 \le i_p \le m$. Then there is another length $q$ cycle, say $C_2$,
such that $|v_2 - i_p| \le m$. Following the same above argument, we conclude that
$C_1 C_2$ is a product of $p+q$ transpositions in $G_m$. By a similar argument as above,
one can conclude that this contradicts  our assumption.
Therefore, the elements of all $C_i$, $1 \le i \le d$, lie in $\{1,\dots,d+1\} \cup \{n-d,\dots,n\}$.
Furthermore, at most two cycles have length greater than $2$.

{\bf Case I} Suppose all these cycles have length $2$. Then by (\ref{uv1}) and (\ref{uv2}),
the case (b.1) follows.

{\bf Case II}  Suppose all cycles $C_1,\dots,C_d$ have length at most $3$.
Since all these cycles are in $L_m$, each of them
can contain at most one of $d+1$ or $n-d$ but not both.
We claim that exactly one of $d+1$ or $n-d$ does not lie in any of cycles.
Suppose not, that is, there are two cycles and each of them contains $d+1$ and $n-d$ respectively.
If both of two these cycles are of length $2$, say $C_1 = (i_1,i_2)$ and $C_2 = (j_1,j_2)$ with
$j_1 < i_1 = d+1 < j_2 = n-d < i_2$.
Then $C_1 C_2 = (i_1,j_2) (j_1,i_1) (j_2,i_2) (i_1,j_2)$
is a product of $4$ transpositions in $G_m$.
By Lemma \ref{2.9}, $\sigma$ is a product of $n+d-3$ transpositions in $G_m$, a contradiction.
Now if one of these two cycles has length $3$ while another has length $2$, say
$C_1 = (i_1,i_2)$ with $i_1 = d+1$ and $i_2 \ge n-d+1$ and $C_2 = (j_1,j_2,j_3)$ with $j_2 = n-d$.
Then by (\ref{uv2}), either $j_1 \le d < j_2 < j_3$ or $j_3 \le d < j_2 < j_1$. Then
$$C_1 C_2 =  \begin{cases}
(j_2,j_3) (i_1,j_2) (j_1,i_1) (j_2,i_2) (i_1,j_2) & \hbox{if } j_1 < j_2 < j_3 \cr
(i_1,j_2) (j_3,i_1) (j_2,i_2) (i_1,j_2) (j_2,i_1)  & \hbox{if } j_3 < j_2 < j_1
\end{cases}$$
In both cases, $C_1C_2$ can be written as a product of $5$ transpositions in $G_m$
and hence  by Lemma \ref{2.9} $\sigma$ is a product of $n+d-3$ transpositions in $G_m$,
a contradiction.
Similar argument can show that it is impossible to have a length $2$ cycle containing $n-d$ while
another length $3$ cycle containing $d+1$.
Finally if there are two length $3$ cycles containing $d+1$ and $n-d$ respectively,
say $C_1 = (i_1,i_2,i_3)$ and $C_2 = (j_1,j_2,j_3)$ with $i_2 = d+1$ and $j_2 = n-d$.
By (\ref{uv1}) and (\ref{uv2}), we may further assume that $\{u_1,v_1\} = \{i_1,i_3\}$
and $\{u_2,v_2\} = \{j_1,j_3\}$. Then
$$C_1 C_2 = \begin{cases}
(i_1,i_2)(i_2,j_2)(i_2,j_1)(i_3,j_3)(i_2,j_2)(j_2,j_3)
& \hbox{if } i_1 < i_2 < i_3 \hbox{ and } j_1 < j_2 < j_3, \cr
(i_1,i_2) (j_1,j_2) (i_2,j_2) (i_2,j_3) (i_3,j_2) (i_2,j_2)
& \hbox{if } i_1 < i_2 < i_3 \hbox{ and } j_3 < j_2 < j_1, \cr
(i_2,j_3) (i_2,j_1) (i_1,j_2) (i_2,j_2) (i_2,i_3) (j_2,j_3)
& \hbox{if } i_3 < i_2 < i_1 \hbox{ and } j_1 < j_2 < j_3, \cr
(j_1,j_2) (i_2,j_2) (i_1,j_2) (i_3,j_3) (i_2,j_2) (i_3,j_3)
& \hbox{if } i_3 < i_2 < i_1 \hbox{ and } j_3 < j_2 < j_1.
\end{cases}$$
In all cases, $C_1C_2$ is a product of $6$ transpositions in $G_m$,
and hence $\sigma$ is a product of $n+d-3$ transpositions in $G_m$,
which contradicts our assumption.
Therefore, we conclude that one and only one of $d+1$ and $n-d$ does not lie in any
of $C_1,\dots,C_d$.
Hence, exactly one of the cycles has length $3$ and all other cycles have length $2$.
Then (b.2) holds if one of the cycles contains $d+1$ and (b.3) holds otherwise.

{\bf Case III} Suppose exactly one cycle has length $4$. Then all other cycles have length $2$.
By (\ref{uv1}) and (\ref{uv2}), the length $4$ cycle must contain at least one of $d+1$ and $n-d$.
Suppose the cycle $C_1 = (i_1,i_2,i_3,i_4)$ contains only one of them, say $n-d$.
Let $i_1 = n-d$. Then there is another length $2$ cycle, say $C_2 = (j_1,j_2)$
with $j_1 = d+1$ and $j_2 \ge n-d+1$. If $|i_4 - i_3| \le m$, then
$(i_3,i_4) C_1 = (i_1,i_2,i_3)$  is a length $3$ cycle containing $n-d$.
Then $(i_3,i_4) C_1$ and $C_2$ satisfy the condition (\ref{rs}).
By Lemma \ref{2.9}, $(i_3,i_4) C_1C_2$ can be written as a product of $5$ transpositions in $G_m$
and so $C_1C_2$ is a product of $6$ transpositions.
Applying Lemma \ref{2.9} again one can conclude that $\sigma$ is a product of $n+d-3$
transpositions in $G_m$, a contradiction. Therefore, $|i_4 - i_3 | > m$.
Similarly, one can show that the other two absolute values $|i_3 - i_2|$ and  $|i_2-i_1|$
are strictly greater than $m$.
Then one must have $i_1,i_3\in \{1,\dots, d\}$ and $i_2 \in \{n-d+1,\dots,n\}$.
Thus, (b.4) follows.  By the similar argument, the result holds if $C_1$ contains only $d+1$
but not $n-d$. Finally, suppose $C_1$ contains both $d+1$ and $n-d$.
If $\{i_1,i_2\} = \{d+1,n-d\}$, then
$(i_3,i_4) C_1 = (i_1,i_2,i_3)$ contains both $d+1$ and $n-d$ and so it is not in $L_m$.
By Lemma \ref{2.9}, it is a product of $2$ transpositions in $G_m$ and
so $C_1$ is the product of $3$ transpositions in $G_m$.
Thus, $\sigma$ is a product of $n+d-2$ transpositions in $G_m$, a contradiction.
So $\{i_1,i_2\}$ contain at most one of $d+1$ and $n-d$.
The same observation holds for $\{i_2,i_3\}$, $\{i_3,i_4\}$ and $\{i_1,i_4\}$.
It follow that either $\{i_1,i_3\} = \{d+1,n-d\}$ or $\{i_2,i_4\} = \{d+1,n-d\}$.
Thus, (b.4) holds.

It remains to show that all the permutations mentioned in (a) and (b) attain the upper bound.
Suppose $n-m$ is even and $\sigma$ has the required form in (a).  Let $C_t=(i_t,j_t)$ for
$t=1,\ldots,d$ and $C_{d+1}=(k_1,\ldots,k_{n-2d})$. Suppose $d({\bf 1},\sigma,m)=s$ and
$\tau_1,\ldots,\tau_s$ are transpositions in $G_m$ such that
$$\tau_1\cdots\tau_s\sigma=1.$$
Assume that $\tau_v\in K_2(\tau_{v+1}\cdots \tau_s\sigma)$ for $v=r_1,\ldots,r_q$ with
$1\le r_1<\cdots<r_q\le s$ and
$\tau_v\in K_1(\tau_{v+1}\cdots \tau_s\sigma)$ for other $v$'s.
Since $C_1,\ldots,C_d\not\in G_m$ and $\min_{1\le t\le d}j_t-\max_{1\le t\le d}i_t=(n-d+1) - d > m$,
one needs at least one term in $\{d+1\ldots,n-d\}$ and   $d$ distinct type two transpositions
to move elements in $C_1,\ldots,C_d$ back to their natural positions. Thus we have $q\ge d$.
Notice that
$$d({\bf 1},\tau_{r_1}\cdots \tau_s\sigma,m)=r_1-1.$$
On the other hand, by Lemma 2.1, the number of disjoint cycles in $\tau_{r_1}\cdots \tau_s\sigma$
is
$$d+1+(s-r_1+1-q)-q=d-2q+s-r_1+2,$$
which implies
$$d({\bf 1},\tau_{r_1}\cdots \tau_s\sigma,n-1)=n-(d-2q+s-r_1+2).$$
It follows that
\begin{eqnarray*}
r_1-1=d({\bf 1},\tau_{r_1}\cdots \tau_s\sigma,m)
\ge d({\bf 1},\tau_{r_1}\cdots \tau_s\sigma,n-1)=n-(d-2q+s-r_1+2)\ge n-(-d+s-r_1+2)
\end{eqnarray*}
which ensures $d({\bf 1},\sigma,m)=s\ge n+d-1$.

When $n-m$ is odd and $\sigma$ has the required form in (b.1), (b.2) or (b.3), one can use the same above argument to deduce that  $\sigma$ attains the upper bound.

Now suppose $\sigma$ has the required form in (b.4).
Denote by $C_t=(i_t,j_t)$ for $t=1,\ldots,d-1$,
$C_d=(i_d,j_d,i_{d+1},j_{d+1})$ and $C_{d+1}=(k_1,\ldots,k_{n-2d-2})$.
Again, suppose $d({\bf 1},\sigma,m)=s$ and $\tau_1,\ldots,\tau_s$ are transpositions in $G_m$
such that
\begin{equation}\label{s11}
\tau_1\cdots\tau_s\sigma=1.
\end{equation}
 Note that $\{i_1,\dots,i_{d+1}\} = \{1,\dots,d+1\}$,
$\{j_1,\dots,j_{d+1}\} = \{n-d,\dots,n\}$. There are at least $d-1$ type two transpositions
in $\tau\equiv\{\tau_1,\ldots,\tau_s\}$. If there are $d$ type two transpositions in $\tau$, then
the same argument  for (a) works for (b.4).

Next we assume there are exactly $d-1$ type two transpositions in $\tau$.
Then $C_1,\ldots,C_{d}$ can be reordered as
$C_{p_1},\ldots,C_{p_v},C_{p_{v+1}},\ldots,C_{p_d}$ with $v\ge 1$
such that there is no transposition  $(u,w)\in \tau$ with
$$u\in V(C_{p_1},\ldots,C_{p_v}), w\in V(C_{p_{v+1}},\ldots,C_{p_d},C_{d+1}),$$
where $V(C_{p_1},\ldots,C_{p_v})$ denotes the set of elements  in cycles $C_{p_1},\ldots,C_{p_v}$.
Since
$$ V(C_{p_1},\ldots,C_{p_v})\cap V(C_{p_{v+1}},\ldots,C_{p_d},C_{d+1})=\emptyset,$$ by (\ref{s11})
we can assume
\begin{eqnarray*}
\tau_1\cdots\tau_u C_{p_1}\cdots C_{p_v}={\bf 1}
\end{eqnarray*}
and
\begin{eqnarray}\label{s12}
\tau_{u+1}\cdots\tau_s C_{p_{v+1}}\cdots C_{p_{d}}C_{d+1}={\bf 1} 
\end{eqnarray}
for some $u$.
Notice that to move elements in $ V(C_{p_1},\ldots,C_{p_v})$ back to their natural positions, we
must  use the transposition $(d+1,n-d)$. Since  $\{d+1,n-d\} \cap \{i_d,j_d,i_{d+1},j_{d+1}\}\ne \emptyset$, we have $C_{d}\in  \{C_{p_1},\ldots,C_{p_v}\}$ and
$ V(C_{p_1},\ldots,C_{p_v})$ contains $2(v+1)$ elements. Moreover, among the transpositions
$ \tau_1,\ldots,\tau_u$, there are at least $v+1$ transpositions with  form $(d+1,n-d)$, $v$
transpositions with form $(i,d+1)$, $1\le i\le d$, and $v$ transpositions with form $(j,n-d)$,
$n-d+1\le j\le n$. Hence
\begin{equation}\label{s13}
u\ge 3v+1.
\end{equation}

 Suppose there are $q$ type two transpositions in $\tau_{u+1},\ldots,\tau_s$.
 By (\ref{s12}),  $q\ge d-v$. Assume the first type two transpositions among them is $\tau_{r}$.
 Denote by $\alpha=\tau_{r}\cdots\tau_sC_{p_{v+1}}\cdots C_{p_{d}}C_{d+1}$.
 Then $$d({\bf 1},\alpha,m)=r-1-u.$$
On the other hand, the number of disjoint cycles in $\alpha$ is
$$d+1-v+(s-r+1-q)-q=d-2q-r+s-v+2.$$
Therefore,
$$d({\bf 1}, \alpha, n-1)=n-2v-2-(d-2q-r+s-v+2).$$
It follows that
$$r-1-u=d({\bf 1},\alpha,m)\ge d({\bf 1}, \alpha, n-1)=n-2v-2-(d-2q-r+s-v+2).$$
Hence,
$$s-u\ge n+2q-d-v-3\ge n+d-3v-3.$$
By (\ref{s13}) we have
$$d({\bf 1},\sigma,m)=s\ge n+d-2.$$
On the other hand, by the former arguments, $\sigma$ can be written as a product of $n+d-1$
transpositions from $G_m$. By  even and odd permutation
rule, we have
$$d({\bf 1}, \sigma,  m)\geq n+d-1. $$
\vskip -.4in \qed

\subsection{The case when $1 < m < (n-1)/2$}

Since $m\ge 2$, we need only to discuss the case when $n\ge 6$.
We are not able to determine the exact value of $\delta(n,m)$
for these cases. Nevertheless,
we have the following upper bounds.

\begin{proposition} \label{2.5}
Let $n, m$ be integers with $n\geq 6$ and $1 < m< \frac{n-1}{2}$. Then
$$\delta(n,m) \le \left\lceil \frac{n-1}{m} \right\rceil + \delta(n-1,m).$$
\end{proposition}

\it Proof. \rm  Let $\sigma\in S_n$ with $\sigma(i)=n$.
Suppose $n-i=qm+r$ with $0 < r \le m$. Take the transpositions
$$(i,i+m),(i+m,i+2m),\ldots,(i+(q-1)m,i+km),(i+qm,n).$$
Thus we move $n$ to the last position by $q+1 = \left\lceil \frac{n-i}{m} \right\rceil$
transpositions and get a new permutation $\sigma_1=[j_1,\ldots,j_{n-1},n]$. Note that
$\tau\equiv[j_1,\ldots,j_{n-1}]$ is a permutation in $S_{n-1}$. We have
\begin{eqnarray*}
d({\bf 1}, \sigma,  m)&\leq&  \left\lceil \frac{n-i}{m}\right\rceil +d({\bf 1}, \sigma_1, m)
= \left\lceil \frac{n-i}{m}\right\rceil+d({\bf 1}, \tau, m) \\
&\leq&  \left\lceil \frac{n-i}{m}\right\rceil+\delta(n-1,m)
\leq  \left\lceil \frac{n-1}{m}\right\rceil+\delta(n-1,m).
\end{eqnarray*}
\vskip -.4in \qed

\begin{proposition} \label{2.6}
Let $n, m$ be integers with $n\geq 6$ and $1 <  m< (n-1)/2$. Then $$\delta(n,m) \le
2 \left\lceil \frac{n-1}{m} \right\rceil-1 + \delta(n-2,m).$$
\end{proposition}

\it Proof. \rm
Let $\sigma\in S_n$ with $\sigma(i)=n$ and $\sigma(j)=1$.
It suffices to verify that we can move 1 and $n$ back to their positions in
$2\left\lceil \frac{n-1}{m} \right\rceil-1$ steps and get a new permutation $[1,i_2,\ldots,i_{n-1},n]$.
Since $[i_2,\ldots,i_{n-1}]$  is a permutation in $S_{n-2}$, we can fix it in
at most $\delta(n-2,m)$ steps. Therefore
$$\delta(n,m)=\max_{\sigma\in S_n}d({\bf 1},\sigma,m)\leq 2\left\lceil \frac{n-1}{m}\right\rceil-1 + \delta(n-2,m).$$

Suppose $n-1=mq+r$ with $q,r$ be integers and $0<r\leq m$.
Since $m < (n-1)/2$, we have $q \ge 2$.
Let $s$ and $t$ be positive integers such that
$$(s-1)m+1 \le i < sm+1 \quad\hbox{and}\quad tm + 1 < j \le (t+1)m+1.$$
Notice that $|t-s| \le q-1$. Suppose $s \ge t$.
Then we can use the following  transpositions to move $1$ and $n$ back to
the first position and the last position, respectively
\begin{multline*}
(i,sm+1),\ (sm+1,(s+1)m+1),\ \ldots,\ ((q-1)m+1,qm+1),\ (qm+1,n),\\
(j, tm+1),\ (tm+1,(t-1)m+1),\ \dots,\ (2m+1,m+1),\ (m+1,1).
\end{multline*}
The number of these transpositions is
$$(q-s+2) +(t+1) = q + 3 + t - s \le q+3  \le 2q + 1 =  2\left\lceil \frac{n-1}{m}\right\rceil-1.$$
Suppose $s < t$.
Then we can use the following  transpositions to move $1$ and $n$ back to
the first position and the last position, respectively
\begin{multline*}
(j,tm+1),\ (i,sm+1),\
(sm+1,(s+1)m+1),\ \ldots,\ ((q-1)m+1,qm+1),\ (qm+1,n),\\
((t-1)m+1,(t-2)m+1),\ \dots,\ (2m+1,m+1),\ (m+1,1).
\end{multline*}
The number of these transpositions is
$$(q-s+3)+(t-1) = q+2 + t-s \le q+2+(q-1) = 2q+1 =  2\left\lceil \frac{n-1}{m}\right\rceil-1.$$
\vskip -.4in \qed

\subsection{Results on $\delta(n,m)$ for small $n$}

By the results in the previous sections and numerical computation,
we have the following table for $\delta(n,m)$.
$$\begin{array}{c|ccccccccccc}
{n\backslash m} & 1&2&3&4&5&6&7&8&9&10&11  \cr
\hline
2&1                &&&&&&&&&& \cr
3&3&2               &&&&&&&&& \cr
4&6&4&3              &&&&&&&& \cr
5&10&5&5&4            &&&&&&& \cr
6&15&[7]&6&6&5         &&&&&& \cr
7&21&[10]&8&7&7&6       &&&&& \cr
8&28&[14]&[10]&9&8&8&7   &&&& \cr
9&36&[16]&[11]&10&10&9&9&8&&& \cr
10&45&[19]&[14]&[12]& 11&11&10&10&9&&\cr
11&55&[23]&[16]&[14]&13&12&12&11&11&10&\cr
12&66&29^*&20^*&17^*&16^*&14&13&13&12&12&11\cr
\end{array}
$$
\centerline{$\hbox{Table 1}$}

Here, the number marked with ``$*$'' are upper bounds for $\delta(n,m)$.
Note that the upper bounds of the numbers in the table are obtained by
Propositions \ref{2.5} and \ref{2.6}. For example,
by Proposition \ref{2.5},
$\delta(12,2)\leq 29, ~\delta(12,3)\leq 20$; by Propositions
\ref{2.5} and also \ref{2.6},
$\delta(12,4)\leq 16$.

The values in square bracket $[\cdot]$
were computed by a Java program written by the third author.
The program uses breadth first search to generate permutations from
$[1,\dots, n] \in S_n$ using elements in $G_m$, and identifies
the minimum number
of steps needed to generate all permutations in $S_n$, and also the
permutations require the maximum number of steps.
One can download the program source code from the link
``https://github.com/sharonli/permutation''.
When the   command ``Permutation'' is executed, one will be asked to input
$n$ and $m$.  One will also be asked whether the program should show all the permutations
generated in each step. If one says no, only the permutations generated in the final step
will be displayed.

Because of memory limitations, this program can handle
the problem up to $n = 11$.  It is easy to modify the program to determine
the diameter of Cayley graphs with vertices connected by other sets of
permutations.
Some numerical results obtained by the program are shown in the Appendix of the paper.

\section*{Acknowledgement}

The study of the problem in the paper
began when C.K. Li was visiting the Hong Kong University of Science
and Technology supported by a Fulbright Fellowship in 2011.
The problem arose in a general education course  ``Mathematics in Daily Life''
conducted by C.K. Li. He acknowledges the support of Fulbright Foundation,
the Hong Kong University of Science and Technology, and some helpful discussion with
colleagues and students at the Hong Kong University of Science and Technology.
The research of Sze was supported by a Hong Kong RGC grant PolyU 502411.
The research of Huang was supported by the NSFC grant 11401197 and a Fundamental Research
Fund  for the Central Universities. This work  began when
Huang was working as a Postdoctoral Fellow at The Hong Kong Polytechnic University.
He thanks the Hong Kong Polytechnic University for its hospitality and support.
The authors are grateful to the referee for helpful suggestions.

\newpage
\section*{Appendix: Some numerical results}
\

\smallskip
\noindent
1. The same 4 permutations attain $\delta(7,2) = 10$ and $\delta(7,3) =  8$. \\
\hspace*{5mm}
$[6, 7 , 4, 5, 2, 3, 1],
 [6, 7, 4, 5, 3, 1, 2],
 [6, 7, 5, 3, 4, 1, 2],
 [7, 5, 6, 3, 4, 1, 2].$

\medskip\noindent
2. The same unique permutation attains $\delta(8,2) = 14$ and $\delta(8,3) = 10$, namely,
$[8,7,6,5,4,3,2,1].$

\medskip\noindent
3. There is a unique permutation attaining $\delta(9,2) = 16$, namely,
$[9,8,7,6,5,4,3,2,1].$

\medskip\noindent
4. There are $770$ permutations attaining $\delta(9,3) = 11$.

\medskip\noindent
5. There are 39 permutations attaining $\delta(10,2) = 19$.\\
\hspace*{5mm}
$[10,9,8,7,6,5,4,3,2,1],
[9,10,8,7,5,6,3,4,1,2],[10,9,8,6,7,4,5,3,2,1],$ \\
\hspace*{5mm}
$[10,9,8,7,6,4,5,2,3,1],
[9,10,8,7,6,5,3,4,2,1],
[10,9,8,7,6,5,2,4,3,1],$\\
\hspace*{5mm}
$[9,10,7,8,5,6,3,4,2,1],[10,9,7,8,5,6,4,3,2,1],
[10,9,8,7,4,6,5,3,2,1],$ \\
\hspace*{5mm}
$[10,8,9,6,7,4,5,2,3,1],
[10,9,6,8,7,5,4,3,2,1],
 [10,9,8,7,5,4,6,3,2,1],$ \\
\hspace*{5mm}
$[10,9,7,8,5,6,3,4,1,2],
[10,9,8,7,6,5,4,1,3,2],
[10,9,8,6,7,5,4,2,3,1],$ \\
\hspace*{5mm}
$ [10,8,9,7,6,4,5,3,2,1],
[10,9,8,7,5,6,3,4,2,1],
[9,10,7,8,6,5,4,3,2,1],$ \\
\hspace*{5mm}
$[9,10,8,7,6,5,4,3,1,2],[10,8,7,9,6,5,4,3,2,1],
[10,9,8,7,6,3,5,4,2,1],$ \\
\hspace*{5mm}
$[9,8,10,7,6,5,4,3,2,1],
[9,10,7,8,5,6,4,3,1,2],
 [10,9,8,7,6,5,4,2,1,3],$ \\
\hspace*{5mm}
$[8,10,9,7,6,5,4,3,2,1],
[10,9,8,7,6,5,3,4,1,2],
[10,9,8,7,5,6,4,3,1,2],$ \\
\hspace*{5mm}
$[10,9,7,8,6,5,3,4,2,1],
[10,8,9,7,6,5,4,2,3,1],
[10,9,7,6,8,5,4,3,2,1],$ \\
\hspace*{5mm}
$[10,8,9,6,7,5,4,3,2,1],
[10,7,9,8,6,5,4,3,2,1],
[10,9,7,8,6,5,4,3,1,2],$ \\
\hspace*{5mm}
$[10,9,8,5,7,6,4,3,2,1],
[10,9,8,7,6,5,3,2,4,1],
[10,9,8,6,5,7,4,3,2,1],$ \\
\hspace*{5mm}
$[9,10,7,8,6,5,3,4,1,2],
[9,10,8,7,5,6,4,3,2,1],
[10,9,8,7,6,4,3,5,2,1].$

\medskip\noindent
6. There are 8 permutations attaining $\delta(10,3) = 14$. \\
\hspace*{5mm}
$[9,10,7,8,6,5,4,3,1,2],
[10,9,7,8,6,5,4,3,2,1],
[10,9,8,7,6,5,3,4,2,1],$ \\
\hspace*{5mm}
$
[9,10,8,7,6,5,3,4,1,2],[10,9,7,8,6,4,5,1,2,3],
[8,9,10,5,6,7,3,4,2,1],$ \\
\hspace*{5mm}
$
[8,9,10,6,7,5,3,4,2,1],
[10,9,7,8,4,5,6,1,2,3].$

\medskip\noindent
7. There are 38 permutations attaining $\delta(10,4) = 12$. \\
\hspace*{5mm}
$[8,9,10,6,7,5,4,1,2,3],
[8,9,10,7,4,5,6,1,2,3],
[10,7,8,9,6,5,4,3,2,1],$ \\
\hspace*{5mm}
$
[9,8,10,6,7,5,4,2,1,3],[8,10,9,7,4,5,6,1,3,2],
[10,8,9,5,7,4,6,2,3,1],$ \\
\hspace*{5mm}
$
[8,10,9,6,4,7,5,1,3,2],
[10,9,8,7,6,4,5,3,2,1],[8,9,10,6,4,7,5,1,2,3],$ \\
\hspace*{5mm}
$
[8,10,9,5,7,4,6,1,3,2],
[8,9,10,5,7,4,6,1,2,3],
[8,10,9,7,6,4,5,1,3,2],$ \\
\hspace*{5mm}
$[10,8,9,6,7,5,4,2,3,1],
[10,8,9,7,4,5,6,2,3,1],
[10,8,9,6,4,7,5,2,3,1],$ \\
\hspace*{5mm}
$
[9,8,10,5,6,7,4,2,1,3],[8,10,9,5,6,7,4,1,3,2],
[9,10,8,6,4,7,5,3,1,2],$ \\
\hspace*{5mm}
$
[9,10,8,6,7,5,4,3,1,2],
[9,8,10,7,6,4,5,2,1,3],[10,9,8,5,6,7,4,3,2,1],$ \\
\hspace*{5mm}
$
[10,9,8,5,7,4,6,3,2,1],
[10,9,8,7,6,5,2,3,4,1],
[9,8,10,6,4,7,5,2,1,3],$ \\
\hspace*{5mm}
$[8,9,10,7,6,4,5,1,2,3],
[8,9,10,5,6,7,4,1,2,3],
[9,8,10,5,7,4,6,2,1,3],$ \\
\hspace*{5mm}
$
[9,10,8,7,4,5,6,3,1,2],[10,8,9,7,6,4,5,2,3,1],
[9,10,8,5,6,7,4,3,1,2],$ \\
\hspace*{5mm}
$
[9,10,8,7,6,4,5,3,1,2],
[8,10,9,6,7,5,4,1,3,2],[9,8,10,7,4,5,6,2,1,3],$ \\
\hspace*{5mm}
$
[10,8,9,5,6,7,4,2,3,1],
[9,10,8,5,7,4,6,3,1,2],
[10,9,8,7,4,5,6,3,2,1],$ \\
\hspace*{5mm}
$[10,9,8,6,4,7,5,3,2,1],
[10,9,8,6,7,5,4,3,2,1].$

\medskip\noindent
8. There are 19 permutations attaining $\delta(11,2) = 23$. \\
\hspace*{5mm}
$[11,10,9,8,7,6,5,4,3,2,1], [11,10,9,8,7,6,5,2,4,3,1],
[11,10,9,7,6,8,5,4,3,2,1],$ \\
\hspace*{5mm}
$[11,10,9,8,7,6,5,4,2,1,3],
[11,10,7,9,8,6,5,4,3,2,1],
[9,11,10,8,7,6,5,4,3,2,1],$ \\
\hspace*{5mm}
$[11,10,9,8,7,6,5,3,2,4,1],
[11,10,9,6,8,7,5,4,3,2,1],
[11,8,10,9,7,6,5,4,3,2,1],$ \\
\hspace*{5mm}
$[11,9,8,10,7,6,5,4,3,2,1],
[11,10,9,8,7,6,5,4,1,3,2],
[11,10,8,7,9,6,5,4,3,2,1],$ \\
\hspace*{5mm}
$[11,10,9,8,7,4,6,5,3,2,1],
[11,10,9,8,7,5,4,6,3,2,1],
[10,9,11,8,7,6,5,4,3,2,1],$ \\
\hspace*{5mm}
$[11,10,9,8,6,5,7,4,3,2,1],
[11,10,9,8,7,6,3,5,4,2,1],
[11,10,9,8,5,7,6,4,3,2,1],$ \\
\hspace*{5mm}
$
[11,10,9,8,7,6,4,3,5,2,1]$.

\medskip\noindent
9. There are 170 permutations attaining $\delta(11,3) = 16$. \\
\hspace*{5mm}
$[8,10,9,11,6,7,5,4,2,3,1],
[10,11,8,9,6,7,5,4,3,2,1],
[11,10,9,5,6,7,8,4,2,3,1],$\\
\hspace*{5mm}
$[11,9,7,8,6,10,2,4,5,3,1],
[11,7,10,8,9,6,5,4,2,3,1],
[11,9,10,8,6,7,3,5,4,1,2],$\\
\hspace*{5mm}
$[9,11,10,8,4,5,6,7,2,3,1],
[11,9,10,8,7,5,6,2,3,4,1],
[11,10,9,8,7,6,5,4,2,3,1],$\\
\hspace*{5mm}
$[11,10,9,8,6,5,7,4,2,3,1],
[11,8,9,10,6,7,5,3,4,2,1],
[11,9,10,7,8,5,6,4,3,2,1],$\\
\hspace*{5mm}
$[11,9,10,8,7,5,6,4,2,3,1],
[11,9,10,8,7,6,4,5,2,3,1],
[8,9,10,11,7,5,6,1,2,3,4],$\\
\hspace*{5mm}
$[11,9,8,10,6,7,5,2,4,3,1],
[11,10,9,8,6,7,2,4,5,3,1],
[11,9,10,8,7,6,5,2,4,3,1],$\\
\hspace*{5mm}
$[11,9,10,8,7,5,6,1,3,2,4],
[11,9,10,8,6,5,7,2,4,3,1],
[11,10,9,8,7,5,6,3,4,1,2],$\\
\hspace*{5mm}
$[11,9,8,10,6,7,5,4,3,2,1],
[11,10,8,9,7,5,6,3,4,2,1],
[11,9,10,8,6,7,2,4,3,5,1],$\\
\hspace*{5mm}
$[11,10,8,9,6,7,4,5,3,2,1],
[11,9,10,5,6,7,8,4,3,2,1],
[11,10,5,8,9,7,6,4,3,2,1],$\\
\hspace*{5mm}
$[11,8,10,9,7,5,6,3,2,4,1],
[11,9,7,8,10,2,6,4,5,3,1],
[10,11,7,8,9,5,6,3,4,1,2],$\\
\hspace*{5mm}
$[11,9,10,8,7,6,5,4,3,2,1],
[11,10,7,8,6,5,9,4,3,2,1],
[9,10,11,8,6,7,5,4,2,3,1],$\\
\hspace*{5mm}
$[11,9,10,8,7,5,6,3,2,1,4],
[11,9,10,8,6,5,7,4,3,2,1],
[11,8,7,9,10,5,6,4,2,3,1],$\\
\hspace*{5mm}
$[11,10,9,8,3,5,6,7,4,2,1],
[11,9,10,7,6,5,8,4,2,3,1],
[11,10,7,8,9,5,6,4,2,3,1],$\\
\hspace*{5mm}
$[11,10,9,8,6,7,5,1,2,3,4],
[11,9,10,8,7,2,6,3,5,4,1],
[11,6,8,10,9,7,5,4,2,3,1],$\\
\hspace*{5mm}
$[11,9,7,8,10,5,6,3,2,4,1],
[11,9,10,8,6,7,4,3,5,1,2],
[11,9,10,8,6,7,5,2,3,4,1],$\\
\hspace*{5mm}
$[11,10,5,8,6,7,9,3,4,2,1],
[11,7,10,8,9,5,6,2,4,3,1],
[11,7,9,8,10,5,6,4,2,3,1],$\\
\hspace*{5mm}
$[11,9,10,8,6,5,4,7,2,3,1],
[11,10,9,8,5,7,6,4,2,3,1],
[11,10,6,8,9,7,5,3,4,2,1],$\\
\hspace*{5mm}
$[9,11,10,8,7,5,6,4,1,3,2],
[11,10,9,8,7,5,6,2,4,3,1],
[11,8,9,10,7,5,6,4,2,3,1],$\\
\hspace*{5mm}
$[11,10,9,8,6,7,3,4,2,5,1],
[10,11,9,8,6,7,5,4,2,3,1],
[11,9,10,8,6,7,5,4,2,3,1],$\\
\hspace*{5mm}
$[8,9,10,11,6,7,5,1,2,3,4],
[11,7,10,8,9,5,6,4,3,2,1],
[11,10,9,8,7,5,6,4,3,2,1],$\\
\hspace*{5mm}
$[11,10,7,8,9,6,5,4,3,2,1],
[11,10,8,9,7,3,6,4,5,2,1],
[10,11,7,9,8,5,6,4,2,3,1],$\\
\hspace*{5mm}
$[11,9,10,8,7,5,6,4,1,2,3],
[11,10,8,9,4,5,6,7,3,2,1],
[11,9,10,8,5,7,6,2,4,3,1],$\\
\hspace*{5mm}
$[11,10,8,9,6,7,5,3,4,2,1],
[11,9,10,8,6,7,3,2,4,5,1],
[11,8,10,9,6,7,5,3,2,4,1],$\\
\hspace*{5mm}
$[11,10,8,9,6,7,3,4,5,2,1],
[8,9,10,11,7,5,6,4,3,2,1],
[11,10,9,8,6,5,3,4,7,2,1],$\\
\hspace*{5mm}
$[11,10,8,5,6,7,9,4,3,2,1],
[11,9,10,8,5,7,6,4,3,2,1],
[11,9,8,10,6,7,5,4,2,1,3],$\\
\hspace*{5mm}
$[10,9,11,8,7,5,6,4,2,1,3],
[11,10,9,8,6,7,4,5,2,3,1],
[11,10,9,7,8,5,6,3,4,2,1],$\\
\hspace*{5mm}
$[11,8,10,9,7,6,5,4,2,3,1],
[11,9,10,5,6,7,8,4,2,1,3],
[11,8,10,9,6,5,7,4,2,3,1],$\\
\hspace*{5mm}
$[10,9,11,8,4,5,6,7,2,3,1],
[11,10,9,8,6,7,5,2,4,3,1],
[9,11,10,8,6,7,5,4,1,3,2],$\\
\hspace*{5mm}
$[11,8,9,10,6,7,5,4,2,3,1],
[11,9,10,8,4,6,5,7,2,3,1],
[11,10,9,8,6,7,2,3,4,5,1],$\\
\hspace*{5mm}
$[11,9,7,8,10,6,5,4,2,3,1],
[11,10,9,5,6,7,8,3,4,2,1],
[11,8,10,9,6,7,2,4,5,3,1],$\\
\hspace*{5mm}
$[11,6,10,8,9,7,3,4,2,5,1],
[11,10,9,8,7,6,5,3,4,2,1],
[11,10,9,8,6,7,5,4,3,2,1],$\\
\hspace*{5mm}
$[11,10,8,9,7,5,6,2,3,4,1],
[10,9,11,8,7,5,6,3,2,4,1],
[11,8,10,9,6,7,5,4,1,3,2],$\\
\hspace*{5mm}
$[11,10,9,8,6,7,3,5,4,2,1],
[11,10,9,8,6,5,7,3,4,2,1],
[11,10,9,8,7,6,3,4,5,2,1],$\\
\hspace*{5mm}
$[11,10,8,9,7,5,3,4,6,2,1],
[11,9,10,8,7,6,5,3,2,4,1],
[11,9,10,8,6,7,4,5,3,2,1],$\\
\hspace*{5mm}
$[11,9,10,8,7,5,6,3,4,2,1],
[11,9,10,5,8,7,6,4,2,3,1],
[11,9,8,10,7,6,5,4,2,3,1],$\\
\hspace*{5mm}
$[11,9,10,8,6,5,7,3,2,4,1],
[11,10,8,9,7,5,6,4,2,3,1],
[11,10,7,9,8,5,6,4,3,2,1],$\\
\hspace*{5mm}
$[11,9,8,10,6,5,7,4,2,3,1],
[11,9,10,5,7,6,8,4,2,3,1],
[11,10,7,8,6,9,5,3,4,2,1],$\\
\hspace*{5mm}
$[9,11,10,8,7,5,6,2,4,3,1],
[10,11,8,7,9,5,6,4,2,3,1],
[11,10,9,8,3,7,6,4,5,2,1],$\\
\hspace*{5mm}
$[11,10,9,8,4,5,6,7,2,3,1],
[11,10,9,8,4,5,6,3,7,2,1],
[10,9,11,8,6,7,5,4,2,1,3],$\\
\hspace*{5mm}
$[11,10,8,9,3,5,6,4,7,2,1],
[11,9,10,8,6,7,2,3,5,4,1],
[11,8,10,9,5,7,6,4,2,3,1],$\\
\hspace*{5mm}
$[11,9,8,10,6,7,2,4,5,3,1],
[11,10,8,9,6,7,2,4,3,5,1],
[11,8,10,9,6,7,3,4,2,5,1],$\\
\hspace*{5mm}
$[11,10,5,9,6,7,8,4,3,2,1],
[9,10,11,8,6,7,5,3,4,2,1],
[11,7,10,8,9,5,6,3,2,4,1],$\\
\hspace*{5mm}
$[11,10,8,9,7,6,5,4,3,2,1],
[11,9,10,8,7,6,2,4,5,3,1],
[11,9,10,8,4,7,6,5,2,3,1],$\\
\hspace*{5mm}
$[11,10,9,8,7,5,6,3,2,4,1],
[11,10,8,9,6,5,7,4,3,2,1],
[11,10,7,8,9,5,6,3,4,2,1],$\\
\hspace*{5mm}
$[11,8,10,9,7,5,6,4,3,2,1],
[11,7,8,10,9,5,6,4,2,3,1],
[11,9,10,5,6,7,8,4,1,3,2],$\\
\hspace*{5mm}
$[11,7,10,8,9,5,3,4,2,6,1],
[11,10,9,8,6,7,3,1,5,2,4],
[11,9,10,7,8,6,5,4,2,3,1],$\\
\hspace*{5mm}
$[11,9,10,8,7,5,6,4,3,1,2],
[8,11,10,9,6,7,5,4,2,3,1],
[11,9,7,8,10,5,6,2,4,3,1],$\\
\hspace*{5mm}
$[11,7,9,8,10,5,6,3,4,2,1],
[11,10,9,8,6,7,4,3,5,2,1],
[11,10,9,8,5,7,6,3,4,2,1],$\\
\hspace*{5mm}
$[11,9,10,8,4,5,6,7,3,2,1],
[11,8,7,9,6,10,5,4,2,3,1],
[10,9,8,11,6,7,5,4,2,3,1],$\\
\hspace*{5mm}
$[11,9,7,8,10,5,6,4,3,2,1],
[10,11,8,9,6,7,3,4,5,1,2],
[10,11,9,8,6,7,3,4,5,2,1],$\\
\hspace*{5mm}
$[11,9,10,8,6,7,5,3,4,2,1],
[11,9,10,8,7,5,6,1,4,3,2],
[11,9,10,8,5,7,6,3,2,4,1],$\\
\hspace*{5mm}
$[11,9,10,8,6,7,3,4,5,2,1],
[11,10,8,9,6,7,5,4,2,3,1],
[8,10,7,11,9,5,6,4,3,2,1],$\\
\hspace*{5mm}
$[11,9,8,10,7,5,6,2,4,3,1],
[11,9,8,10,5,7,6,4,2,3,1],
[11,9,10,8,7,5,3,2,4,6,1],$\\
\hspace*{5mm}
$[11,9,8,10,6,7,3,4,2,5,1],
[11,9,10,8,7,6,3,4,2,5,1],
[11,10,8,9,7,5,6,4,1,2,3],$\\
\hspace*{5mm}
$[11,10,7,8,9,5,6,4,3,1,2],
[11,9,8,10,7,5,6,4,3,2,1],
[11,10,9,7,8,5,6,4,2,3,1],$\\
\hspace*{5mm}
$[11,10,8,7,9,5,6,4,3,2,1],
[11,7,8,9,10,5,6,4,3,2,1],
[11,10,9,8,6,7,5,3,2,4,1],$\\
\hspace*{5mm}
$[11,10,8,9,5,7,6,4,3,2,1],
[11,8,10,9,6,7,5,4,3,2,1].$

\medskip\noindent
10. There are 2 permutations attaining $\delta(11,4) = 14$. \\
\hspace*{5mm}
$[10,11,8,9,7,5,6,3,4,1,2],
[10,11,8,9,6,7,5,3,4,1,2].$

\begin{thebibliography}{www}

\bibitem{JD} J.  D$\acute{\rm e}$nes,
The representation of a permutation as the product of a
minimal number of transpositions, and its connection with the theory of graphs.
Magyar Tud. Akad. Mat. Kutat Int. K$\ddot{\rm o}$zl. 4 (1959), 63-71.

\bibitem{E}  J.A.  Eidswick,
Short factorizations of permutations into transpositions.
Discrete Math. 73 (1989), 239-243.

\bibitem{Even}
S. Even,
Graph algorithms. Second edition.
Cambridge University Press, Cambridge, 2012.


\bibitem{VF}
V. F$\acute{\rm e}$ray,
Partial Jucys-Murphy elements and star factorizations.
European J. Combin. 33 (2012),  189-198.

\bibitem{Feng}
X. Feng, B. Chitturi, and H. Sudborough,
Sorting circular permutations by bounded transpositions.
Adv Exp Med Biol. 680 (2010), 725-36.


\bibitem{Fine}
J.T. Fineman and E. Robinson,
Fundamental graph algorithms. Graph algorithms in the language of linear algebra, 4558,
Software Environ. Tools, 22, SIAM, Philadelphia, PA, 2011.

\bibitem{GF}
D.A. Gewurz, F. Merola,
On factorizations of cyclic permutations into transpositions.
Ars Combin. 95 (2010), 397-403.

\bibitem{IR}
J. Irving, A. Rattan,
Minimal factorizations of permutations into star transpositions.
Discrete Math. 309 (2009), 1435-1442.

\bibitem{Knuth}
D. E. Knuth,
The Art of Computer Programming, Volume 3: Sorting and Searching,
Addison-Wesley, Mass.-London-Don Mills, Ont., 1973.


\bibitem{L}
O. P. Lossers,
Solution to Problem E3058,
American Mathematical Monthly 93 (1986), 820-821.


\bibitem{M}
G. Mackiw,
Permutations as products of transpositions.
Amer. Math. Monthly 102 (1995), 438-440.


\bibitem{N}
D. Neuenschwander,
On the representation of permutations as products of transpositions.
Elem. Math. 56 (2001), 1-3.

\bibitem{P}P.
A. Pevzner,
Computational Molecular Biology: An Algorithmic Approach. The MIT Press, Cambridge, 2000.


\bibitem{Yue}
F. Yue, M. Zhang, J. Tang,
Phylogenetic reconstruction from transpositions.
MC Genomics. 2008 Sep 16; 9 Suppl 2:S15.

\bibitem{xiao}
W. Xiao,
Some results on diameters of Cayley graphs.
Discrete Appl. Math. 154 (2006), 1640-1644.

\iffalse\bibitem{iPhone} iPhone applications: Ghost Leg,\\
https://itunes.apple.com/us/app/real-sadari/id352275896.

\bibitem{A} Wikipedia: Ghost Leg, http://en.wikipedia.org/wiki/Ghost\underbar{~}\,Leg.\fi

\end{thebibliography}
\end{document}